\documentclass[AMA,Times1COL]{WileyNJDv5} 
\usepackage{amssymb}
\usepackage{graphicx}
\usepackage{epstopdf} 
\usepackage{hyperref}    
\usepackage{subcaption}
\usepackage{amssymb,latexsym,amsmath,amsthm}
\usepackage{xcolor}
\usepackage{enumerate}
\usepackage{booktabs}
\usepackage{MnSymbol,bbding,pifont}

\usepackage{tikz, pgfplots, mathtools}
\pgfplotsset{compat=1.16}
\theoremstyle{definition}
\newtheorem{example}{Example}[section]

\def\R{\mathbb{R}}
\def\C{\mathbb{C}}
\def\N{\mathbb{N}}

\articletype{Research Article}%

\received{Date Month Year}
\revised{Date Month Year}
\accepted{Date Month Year}
\journal{Numerical linear algebra with applications}
\volume{00}
\copyyear{2024}
\startpage{1}

\begin{document}

\title{Finding the nearest $\Omega$-stable pencil with Riemannian optimization}

\author[1]{Vanni Noferini}

\author[1]{Lauri Nyman}

%\author[2]{Federico Poloni}

\authormark{NOFERINI AND NYMAN}
\titlemark{Nearest $\Omega$-stable pencil}

\address{\orgdiv{Department of Mathematics and Systems Analysis}, \orgname{Aalto University}, \orgaddress{\state{Espoo}, \country{Finland}}}

% \address[2]{\orgdiv{Department Name}, \orgname{Institution Name}, \orgaddress{\state{State Name}, \country{Country Name}}}

\corres{Lauri Nyman, Otakaari 1, P.O. Box 11100, FI-00076, Aalto, Finland. \email{lauri.s.nyman@aalto.fi}}

% \presentaddress{This is sample for present address text this is sample for present address text.}

% \fundingInfo{Text}
% \JELinfo{ejlje}

\abstract[Abstract]{This paper considers the problem of finding the nearest $\Omega$-stable pencil to a given square pencil $A+xB \in \C^{n \times n}$, where a pencil is called $\Omega$-stable if it is regular and all of its eigenvalues belong to the closed set $\Omega$. We propose a new method, based on the Schur form of a matrix pair and Riemannian optimization over the manifold $U(n) \times U(n)$, that is, the Cartesian product of the unitary group with itself. While the developed theory holds for any closed set $\Omega$, we focus on two cases that are the most common in applications: Hurwitz stability and Schur stability. For these cases, we develop publicly available efficient implementations. Numerical experiments show that the resulting algorithm outperforms existing methods. 
}


% \jnlcitation{\cname{%
% \author{Taylor M.},
% \author{Lauritzen P},
% \author{Erath C}, and
% \author{Mittal R}}.
% \ctitle{On simplifying ‘incremental remap’-based transport schemes.} \cjournal{\it J Comput Phys.} \cvol{2021;00(00):1--18}.}

\maketitle

% \renewcommand\thefootnote{}
% \footnotetext{\textbf{Abbreviations:} ANA, anti-nuclear antibodies; APC, antigen-presenting cells; IRF, interferon regulatory factor.}

% \renewcommand\thefootnote{\fnsymbol{footnote}}
% \setcounter{footnote}{1}

\section{Introduction}

In control theory, and especially in the context of studying the asymptotic stability of differential-algebraic equations (DAEs) \cite{dae}, one is motivated to compute the distance from a pencil, i.e., a matrix pair, to the nearest one which is regular and has all its eigenvalues in a given set $\Omega$. In particular, the case where $\Omega$ is the left half-plane is linked to the asymptotic stability of the solution of systems of linear DAEs with constant coefficients \cite{BN,du,gms}; the problem is also complementary to that of the nearest singular pencil \cite{singpencil}. The special case of a single matrix, or equivalently a pencil whose leading coefficient is the identity matrix, is well studied and corresponds to the problem of computing the nearest $\Omega$-stable matrix, i.e., the nearest matrix whose eigenvalues all lie in a closed set $\Omega$. The nearest $\Omega$-stable matrix problem could be seen as a structured variant of the quest for the nearest pencil with eigenvalues in $\Omega$ (where only one matrix coefficient is allowed to be perturbed). Various numerical methods for the nearest $\Omega$-stable matrix problem have been developed \cite{cgs, curtis,GS17,GL17,np,ONV}, especially for $\Omega$ being either the left half plane or the unit disc, but occasionally also for different choices of the set $\Omega$ \cite{cgs,np}. Still in the single-matrix setting, the dual problem of finding the nearest unstable matrix has also received much attention \cite{BM19,Bye88,GM15,kms,mengi}. However, the general problem of computing the nearest $\Omega$-stable pencil is much less studied, and (up to our knowledge) the algorithm by Gillis, Mehrmann and Sharma \cite{gms} is currently the only one available numerical method. The main goal of the present paper is to fill this gap, by extending to the case of pencils the algorithm by Noferini and Poloni \cite{np} to find the nearest $\Omega$-stable matrix. The method in \cite{np}, and hence the method in our paper, is based on a stiffen-and-optimize approach where the original problem is recast as a minimization problem onto a Riemannian manifold. In this paper the manifold in question is $U(n) \times U(n)$, that is the Cartesian product of $U(n)$ (the set of $n \times n$ unitary matrices) with itself. This is contrast to \cite{np}, where the underlying manifold was either $U(n)$ or $O(n)$, and akin to \cite{singpencil} where working on $U(n) \times U(n)$ (or $O(n) \times O(n)$ for the real case) allowed Dopico, Noferini and Nyman to tackle the problem of finding the nearest singular pencil.

As discussed in \cite{gms}, in the general case of pencils the precise definition of stability is not fully consistent in the literature, and for this reason it is useful to clarify from the beginning which problem exactly we are going to tackle. One of the most recent papers on the subject is \cite{gms}, and there the authors call a matrix pair \emph{stable} if the corresponding pencil is regular, of index at most one (i.e, there are no Jordan chains longer than one for infinite eigenvalues), and all its finite eigenvalues are in the closed left half plane, and those that lie on the imaginary axis are semisimple; and they call it \emph{asymptotically stable} if the corresponding pencil is regular, of index at most one, and all its finite eigenvalues are in the open left half plane. In this paper, given a set $\Omega \subseteq \C \cup \{\infty\}$, we call a pencil $\Omega$-stable if it is regular and all its eigenvalues belong to $\Omega$.  Although this definition does not require the set $\Omega$ to be closed, in the present manuscript we are interested in the computation of the nearest $\Omega$-stable pencil of a given one as well as the projection onto the set of $\Omega$-stable pencils; since the distance to a set coincides with the distance to the closure of such set, and the projection is only guaranteed to exist in the closure, mathematically one would anyway need, as a first step, to relax the requirement so that eigenvalues are also allowed to belong to the closure of $\Omega$. As such, we make with no loss of generality the assumption that $\Omega$ is closed. This means that, for the case of Hurwitz stability, we consider $\Omega$-stable pencils with $\Omega=\{ z \in \C : \Re(z) \leq 0 \} \cup \{ \infty \}$; note that this is a superset of both asymptotically stable matrix pairs and stable matrix pairs in the sense of \cite{gms}. We also consider the case of Schur stable pencils, that is, $\Omega$-stable pencils with $\Omega=\{ z \in \C : |z| \leq 1 \}$; just as Hurwitz stability is associated with the behaviour of continuous-time dynamical systems, Schur stability is linked to discrete-time dynamical systems. For Schur stability, we are not aware of any previous algorithm for the pencil case, but only for the matrix case \cite{cgs,np}.

The paper is structured as follows. In Section \ref{sec:problem}, we expand the discussion above and define the problem in a rigorous way. Section \ref{sec:reformulation} follows the lead of \cite{np} and explains how and why the nearest $\Omega$-stable pencil problem is equivalent to an optimization task on a certain Riemannian manifold; later on, in Section \ref{sec:optimization}, we show how to numerically tackle this Riemannian optimization problem for two important special choices of $\Omega$ (corresponding to Hurwitz stable pencils and Schur stable pencils). Before being able to do that, in Section \ref{sec:projformula}, we must solve the technical subproblem of projecting a pair of complex numbers given certain constraints on their ratio. We present numerical experiments in Section \ref{sec:numexp}, for both Hurwitz and Schur stable pencils. In the experiments, we study the problem in general, including a statistical observation that our computed minimizers have nontrivial Jordan chains with positive (and in fact, heuristically, quite high) probability; this generalizes an analogous observation made in \cite{np} for the nearest stable matrix problem. Moreover, for Hurwitz stability, we compare our approach with the state-of-the art algorithm of \cite{gms}, and the results turn out to be significantly in favour of our newly proposed method. Finally, we draw some conclusions in Section \ref{sec:conclusions}.

\section{The nearest stable pencil problem}\label{sec:problem}

A pencil over $\C$, the field of complex numbers, is a polynomial of degree at most $1$ whose coefficients are complex matrices of the same size. In this paper, we consider square pencils, and define $\C[x]_1^{n\times n}$ to be the vector space of $n \times n$ pencils over $\C$, that is, the space  $\C[x]_1^{n \times n}:=\{A+x B$ $|$ $A,B \in \C^{n \times n}\}$. We call a square pencil $A + x B$ \textit{regular} if $\det(A + x B) \not\equiv 0$. Otherwise, we call the pencil \textit{singular}. For a regular pencil, we call a value $x_0 \in \C$ a \textit{finite eigenvalue} of $A + x B$ if, at $x_0$, the pencil evaluates to a singular matrix, that is, $\det(A + x_0 B) = 0$. The \textit{algebraic multiplicity} of $x_0$ as a finite eigenvalue of a regular pencil $A + x B$ is the multiplicity of $x_0$ as a root of $\det (A + xB)$. Below we give a more general definition of a finite eigenvalue (see, e.g., \cite[Section 2]{NV}) that includes singular pencils and reduces to the simpler characterization above in the regular case.
\begin{definition}\label{def:evals}
    The number $x_0 \in \C$ is a finite eigenvalue of the pencil $A+xB$ if
    \[ \mathrm{rank}_\C (A+x_0 B) < \mathrm{rank}_{\C(x)} (A + x B) = \sup_{x \in \C} \mathrm{rank}_\C (A + x B). \]
\end{definition}
The pencil $A + xB$ has an \textit{infinite eigenvalue} if its reversal $B + x A$ has zero as its eigenvalue; in the regular case, equivalently, infinity is an eigenvalue of the regular pencil $A+xB$ precisely when the matrix $B$ is singular. The algebraic multiplicity of infinity as an eigenvalue is the algebraic multiplicity of zero as an eigenvalue of its reversal.

Clearly, $\C[x]_1^{n \times n} \cong \C^{2n^2}$, and therefore $\C[x]_1^{n \times n}$ can be equipped with the Euclidean norm $\| A + x B \|_F$ whose square is defined as follows.
\begin{equation}\label{eq:pencilnorm}
    \| A + x B \|_F^2 := \sum_{i,j=1}^n (|A_{ij}|^2+|B_{ij}|^2) = \left\| \begin{bmatrix}
        A \\
        B
    \end{bmatrix} \right\|_F^2 = \left\| \begin{bmatrix}
        A & B
    \end{bmatrix}  \right\|_F^2,
\end{equation}
where $\|\cdot \|_F$ is the Frobenius matrix norm and, for $a \in \C$, $|a|$ is the modulus of $a$. The norm \eqref{eq:pencilnorm} induces a distance
\begin{equation}\label{eq:pencildist}
    d(A+x B,S+x T):=\|(A-S)+x(B-T) \|_F.
\end{equation}
Throughout this paper, we always refer to \eqref{eq:pencildist} when speaking of the distance between two pencils.

Let $\overline{\mathbb{C}}:=\C \cup \{\infty\}$ denote the topological closure of $\mathbb{C}$ (i.e., up to a stereographic projection, the Riemann sphere). Let $\Omega$ be a non-empty closed subset of $\overline{\mathbb{C}}$, and define
\begin{align}\label{eq:stableset}
    S(\Omega,n) := \{ A + x B \in \C[x]_1^{n \times n} : \Lambda( A + x B) \subseteq \Omega, \ \mbox{det}( A + x B) \not\equiv 0 \} \subseteq \C[x]_1^{n \times n},
\end{align}
the set of regular $n \times n$ matrix pencils whose eigenvalues all belong to $\Omega$. 
Here, $\Lambda ( A + x B)$ denotes the spectrum of the square pencil $ A + x B$. We record this definition formally below.
\begin{definition}
    Given a closed set $\Omega \subseteq \overline{\C}$, an $n \times n$ pencil $A + x B$ is called $\Omega$-stable if $A + x B \in S(\Omega,n)$ as defined in \eqref{eq:stableset}.
\end{definition}

Observe that $S(\Omega,n)$ is not empty. Indeed, since $\Omega$ is non-empty by construction, there exists an element $\omega \in \Omega$ and therefore $M \in S(\Omega,n)$ where we define either $M:=-\omega I + x I$ (if $\omega \in \C$) or $M:=I+x 0$ (if $\omega=\infty$). Furthermore, $S(\Omega,n)$ is not closed, because the sequence $M_k=\frac{M}{k} \rightarrow 0+x0$ as $k \rightarrow \infty$, but $0+x0 \notin S(\Omega,n)$ (not being a regular pencil). This fact is slightly problematic and prompts a brief digression in analysis. 

Indeed, recall that, given an arbitrary metric space $(\mathcal{M},d)$, the distance of a point $x \in \mathcal{M}$ to a set $\mathcal{A} \subseteq \mathcal{M}$ is defined as $d(x,\mathcal{A}):=\inf_{y \in \mathcal{A}} d(x,y)$. If $\mathcal{A}$ is not closed, the infimum in the definition may not be achieved at any point of $\mathcal{A}$; therefore, we can surely study the distance from $x$ to $\mathcal{A}$ but it may well be pointless to ask for a point in $\mathcal{A}$ nearest to $x$. Fortunately, denoting by $\overline{\mathcal{A}}$ the closure of $\mathcal{A}$, the latter shortcoming is easily cured by noting that, by definition of closure, there exists $y_0 \in \overline{\mathcal{A}}$ such that $d(x,\mathcal{A})=d(x,\overline{\mathcal{A}})=d(x,y_0)$\footnote{In fact, a reader interested in a simple exercise may want to prove that $\overline{\mathcal{A}}=\{x \in \mathcal{M} : d(x,\mathcal{A})=0\}$.}. In other words, \emph{closed} sets always possess not only distances to an arbitrary point, but also (at least) one element where such distance is achieved. 

Let us now go back to our original problem, that corresponds to $\mathcal{M}=\C[x]_1^{n \times n}$ and $\mathcal{A}=S(\Omega,n)$. While the distance of a given pencil $A+xB$ to $S(\Omega,n)$ is certainly well defined, our observations lead us to the conclusion that finding the nearest pencil in $S(\Omega,n)$ to $A+xB$ may not be a well posed problem. Hence, we consider instead the distance of $A+xB$ to the closure of $S(\Omega,n)$, denoted by $\overline{S(\Omega,n)}$. This mathematical sophistication does not change the sought distance, while bringing the advantage that a nearest pencil to the input is now guaranteed to exist in $\overline{S(\Omega,n)}$, and so we may have some hope to compute it. The same argument also justifies our choice to assume, with no loss of generality, that $\Omega$ is closed. Indeed, otherwise it is not difficult to prove that $\overline{S(\Omega,n)}=\overline{S(\overline{\Omega},n)}$ (to show this fact, one can slightly modify the proof of Lemma \ref{lem:samedistance}). In summary, the problem of our interest can now be stated as follows.

\begin{problem}\label{problem}
Given a pencil $A+x B \in \C[x]_1^{n \times n}$, compute both the minimum and a minimizer of the distance function $d(A+x B,S+x T)$ \eqref{eq:pencildist} amongst the closure of all $\Omega$-stable pencils, that is, $S+x T \in \overline{S(\Omega,n)}$.
\end{problem}

As it is often more convenient to deal with the \textit{square} of the Euclidean norm, in the following we equivalently express Problem \ref{problem} as the minimization of the square distance:
\begin{equation}\label{eq:riproblem}
    \min_{S+x T \in \overline{S(\Omega,n)}} [d(A+x B,S+x T)]^2.
\end{equation}

   We conclude this section with some further comments on the role of the closure of the set $S(\Omega,n)$ in Problem \ref{problem}. In abstract, it is clear that once one finds a minimizer of the distance $y_0 \in \overline{\mathcal{A}}$, then for all $\epsilon>0$ there exist $y \in \mathcal{A}$ such that $d(y,y_0)<\epsilon$; that is, in plain language, there exist elements of $\mathcal{A}$ that are arbitrarily close to $y_0$. In practice, in Problem \ref{problem}, one is sometimes interested in (approximately) computing the eigenvalues of the pencil $S+xT$ that minimizes the distance to the input. This may be, for example, to verify that these eigenvalues are indeed in $\Omega$, or for other purposes such as further studying the stability of DAEs. There are two possibilities: either $S+xT \in S(\Omega,n)$ or not. In the former case, $S+xT$ will have $n$ eigenvalues in $\Omega$, as required; but in the latter case, since $\Omega$ is closed, then $S+xT$ must be singular and may not have any eigenvalue at all! (For more details, see the discussion in Section \ref{sec:reformulation}, and in particular Proposition \ref{prop:closuresing}, and Example \ref{exnew}.) We illustrate in the simple Example \ref{ex0} below that this situation is far from being purely theoretical and can indeed happen in practice, even  for apparently harmless instances of the problem.

\begin{example}\label{ex0}
Let $n=1$ and $\Omega=\{0\}$. Then,
\[ S(\Omega,1)=\{0+x z : 0 \neq z \in \C \}. \]
Consider the input pencil $1+x0$; its distance to an arbitrary element of $S(\Omega,1)$, parametrized by the nonzero complex number $z$, is $\sqrt{1+|z|^2}$. Thus, it is easy to see that the pencil nearest to the input in $\overline{S(\Omega,1)}$ is $0+x0 \notin S(\Omega,1)$; note also that  $0+x0$ is not a regular pencil, even if it is a limit point of pencils in $S(\Omega,1)$, e.g., $0+\frac{x}{k}$. Note also that, as $k \rightarrow \infty$, the distance of this family of pencils from $1+x0$ tends to $1$. 
\end{example}

When the distance minimizer $S+xT \notin S(\Omega,n)$, it is plausible that one may be interested in a slightly perturbed pencil that does belong to $S(\Omega,n)$; for example, once again this may happen so that eigenvalues can be computed and verified to be in $\Omega$. In Section \ref{sec:reformulation} we will see that, concretely, our method provides a generalized Schur form of the minimizer $S+xT$, say, $S+xT=Q(S_0+xT_0)Z$. In such form, $S_0+xT_0$ is upper triangular, and it is very easy to obtain a perturbed pencil that belongs to $S(\Omega,n)$. Indeed, it suffices to perturb precisely the zero diagonal elements of $S_0+xT_0$, in such a way to create an eigenvalue belonging to $\Omega$.

\section{Reformulation of the problem}\label{sec:reformulation}

In this section, our aim is to extend the ideas first developed by Noferini and Poloni\cite{np} for the problem of computing the nearest $\Omega$-stable matrix to a given one to Problem \ref{problem}, which is the pencil counterpart of the same problem. In their work, Noferini and Poloni made two simple yet synergistic observations: (i) each matrix can be brought to upper triangular form via the Schur form, and (ii) the eigenstructure of a triangular matrix can be easily characterized. In what follows, we will see that analogous statements are true for matrix pencils as well.

% \textcolor{red}{ATM we do not use $SU(n)$ anywhere.}
% Let $U(n)$ denote the unitary group of degree $n$, and let $SU(n)$ denote the corresponding special unitary group, that is, the set of $n \times n$ unitary matrices with determinant $1$. The following lemma generalizes the observation (i) to the pencil setting.

Let $U(n)$ denote the unitary group of degree $n$, that is, the set of $n \times n$ unitary matrices. The following lemma generalizes the observation (i) to the pencil setting.
\begin{lemma}[Stewart \cite{Stewart}] \label{lemma:schur}
    For any pair $A,B \in \C^{n \times n}$ there exist $Q,Z \in U(n)$ such that $QAZ$ and $QBZ$ are both upper triangular.
\end{lemma}
The upper triangular pencil $QAZ + x QBZ$ given by the matrices in Lemma \ref{lemma:schur} is called a generalized Schur form of the pencil $A + x B$. The generalized Schur form is an important tool for eigenvalue computation; indeed, Definition \ref{def:evals} implies that the pencils $A+xB$ and $QAZ+xQBZ$ have the same eigenvalues. In addition, it is also clear that $A+xB$ is regular if and only if $QAZ+xQBZ$ is.
%It is easy to see that each pencil has infinitely many generalized Schur forms. Indeed, a Schur form can be multiplied by arbitrary unitary diagonal matrices on both sides to yield another Schur form; moreover, diagonal elements can be swapped  (see Corollary \ref{cora2} for the regular case).

Clearly, as noted for example in \cite{singpencil}, an upper triangular pencil is singular if and only if at least one of its diagonal entries is equal to $0+x0$. Next, we observe that it is easy to characterize the eigenstructure of regular triangular pencils. This generalizes the observation (ii). To state Lemma \ref{lemma:eigen}, we adopt the convention that if $0 \neq a \in \C$ then $a/0=\infty$.
\begin{lemma} \label{lemma:eigen}
    The multiset of the finite and infinite eigenvalues (counted with algebraic multiplicity) of a regular upper triangular square pencil $A + x B \in \C[x]_1^{n \times n}$ is $\{-a_{ii}/b_{ii} \}_{i=1,\dots,n}$.
\end{lemma}
\begin{proof}
    Note first that, for all $i$, $a_{ii}$ and $b_{ii}$ cannot be both zero due to the regularity assumption, and hence $-a_{ii}/b_{ii}$ is well-defined and equal to either a complex number or $\infty$. At this point, it suffices to observe that the determinant of a triangular pencil is the product of the diagonal elements of the pencil.
\end{proof}
\begin{remark}
    In Lemma \ref{lemma:eigen}, the eigenvalue $-a_{ii}/b_{ii}$ associated with the $i$th diagonal element $a_{ii} + x b_{ii}$ coincides with its unique root. When $b_{ii} = 0$, we say that $a_{ii} + x b_{ii}$ has a \textit{root at infinity}, which corresponds to an infinite eigenvalue of the pencil. With this terminology, the eigenvalues of a regular upper triangular square pencil $A + x B$ are the roots of its diagonal elements.
\end{remark}

Let $\Omega \subseteq \overline{\mathbb{C}}$.
%and define $\C[x]_1^{\Omega} := \{a + x b \in \C[x]_1 : -\frac{a}{b} \in \Omega \}$, that is, $\C[x]_1^{\Omega}$ is the set of linear polynomials with a root in $\Omega$.
%{\color{blue}VN: Is this not just $S(\Omega,1)$? If so I would avoid to introduce two symbols for the same thing.}
We define $p_{\Omega}$ to be the projection from $\C[x]_1$ to the set $\overline{S(\Omega,1)}$, that is,
$$p_{\Omega}: \C[x]_1 \rightarrow \overline{S(\Omega,1)}, \quad p_{\Omega}(a+x b) = \arg \min_{c+x d \in \overline{S(\Omega,1)}} \|a-c\|^2 + \|b-d\|^2. $$
Note that this projection is not necessarily uniquely defined. The set of points at which the projection is not unique is called the \textit{medial axis}. This has consequences for differentiability, as we will see in Lemma \ref{lemma:gradient}. Henceforth, whenever we refer to $p_\Omega$, we tacitly assume that (for those input on which it is not uniquely defined) one possible output has been fixed\footnote{For a general $\Omega$, this step may require the axiom of choice, but we emphasize that this is not needed for common choices of $\Omega$ such as the left half plane or the unit disc.}.

The following proposition shows how we can use Lemma~\ref{lemma:eigen} to find the nearest $\Omega$-stable \textit{upper triangular} pencil.

\begin{proposition} \label{prop:nearest_tri}
    Let $A+x B \in \C[x]_1^{n \times n}$. Let $S_{tri}(\Omega,n)$ be the set of $\Omega$-stable upper triangular pencils, and let $\overline{S_{tri}(\Omega,n)}$ denote its closure. A pencil in $\overline{S_{tri}(\Omega,n)}$ nearest to $A+x B$ is given by $\mathcal{T}(A+ x B)$, where
    \[ \mathcal{T}(A+ x B)_{ij}=\begin{cases}A_{ij} + x B_{ij} \            & \mathrm{if} \ i<j ; \\
             p_{\Omega}(A_{ii}+x B_{ii}) \  & \mathrm{if} \ i=j   \\
             0 \                                  & \mathrm{otherwise}.\end{cases} \]
    % \qquad \mathcal{T}(B)_{ij}=\begin{cases}B_{ij} \ &\mathrm{if} \ i<j ;\\
    % p_{\Omega}(A_{ii}+x B_{ii}) \ &\mathrm{if} \ i=j \\
    % 0 \ &\mathrm{otherwise}. \end{cases}\]
    % In particular, the squared distance to $A+x B$ from $\mathcal{T}(A)+x \mathcal{T}(B)$ is 
    % \begin{equation}\label{eq:objectivefunction} \mathcal{F}(A+x B)= \sum_{i>j} (|A_{ij}|^2+|B_{ij}|^2) + \min_{1 \leq i \leq n} (|A_{ii}|^2+|B_{ii}|^2 ).\end{equation}
\end{proposition}
\begin{proof}
    For upper triangularity, it is necessary that the strictly lower triangular part of the minimizer is zero. An upper triangular pencil is in $\overline{S_{tri}(\Omega,n)}$ if and only if all the diagonal elements belong to $\overline{S(\Omega,1)}$. The projection $p_{\Omega}$ minimizes the distance for each diagonal entry separately. The strictly upper triangular part is left untouched to minimize the total distance.
\end{proof}

%Fix now $A+x B \in \C[x]_1^{n \times n}$. Proposition~\ref{prop:nearest_tri} can be used, with the help of Lemmas~\ref{lemma:schur} and \ref{lemma:singular}, to characterize the singular pencil $S+x T \in \C[x]_1^{n \times n}$ nearest to $A+x B$.  To this goal, recall that $\mathcal{S}_n \subset \C[x]_1^{n \times n}$ is the subset of singular pencils and denote by $\mathcal{T}_n \subset \C[x]_1^{n \times n}$ the subset of singular upper triangular pencils. For any pair of unitary matrices $(Q,Z) \in U(n) \times U(n)$ define the function
%\begin{equation}\label{eq:newobjfun}
% f(Q,Z) := \mathcal{F}(QAZ+x QBZ) = [d(QAZ+x QBZ,\mathcal{T}(QAZ)+x \mathcal{T}(QBZ))]^2  
% \end{equation}
% where $\mathcal{F}$ and $\mathcal{T}$ are as in Proposition~\ref{prop:nearest_tri} and $d$ is as in \eqref{eq:pencildist}. 

Next, we will see how we can approach Problem \ref{problem} by utilizing Lemma~\ref{lemma:schur} and Proposition~\ref{prop:nearest_tri}. To this end, let $A+x B \in \C[x]_1^{n \times n}$ and let $\mathcal{T}$ be as in Proposition~\ref{prop:nearest_tri}. It can be shown that the squared distance from $A+x B$ to the set $\overline{S(\Omega,n)}$ is equal to
\begin{equation}\label{eq:minnewobjfun}
    \min_{(Q,Z)\in U(n) \times U(n)} f(Q,Z),
\end{equation}
where
\begin{equation}\label{eq:newobjfun}
    f(Q,Z) = \| QAZ+x QBZ - \mathcal{T}(QAZ+x QBZ) \|^2_F.
\end{equation}

However, with respect to the analogous developments in the nearest $\Omega$-stable matrix problem \cite{np}, there is a technical complication in the case of pencils. In essence, the reason of the extra difficulty is that the pencil $\mathcal{T}(Ax+B)$ may not be regular (even when $Ax+B$ is regular), and indeed this happens precisely when one of the projections of the diagonal elements of $Ax+B$ is equal to $0+x0$.
Hence, before we can state more formally in Theorem \ref{thm:obj_function} the result that we anticipated above, we need the following lemma.

\begin{lemma} \label{lemma:closures}
    Let $S_{tri}(\Omega,n)$ be the set of $\Omega$-stable upper triangular pencils, and let $\overline{S_{tri}(\Omega,n)}$ denote its closure. It holds that $$\overline{S(\Omega,n)} = \bigcup_{Q,Z \in U(n)} Q \ \overline{S_{tri}(\Omega,n)} \ Z.$$
\end{lemma}
\begin{proof}
For notational simplicity, denote \[\mathcal{S}:=\bigcup_{Q,Z\in U(n)} Q \ \overline{S_{tri}(\Omega,n)} \ Z;\] we aim to show that $\overline{S(\Omega,n)}=\mathcal{S}$.
    The inclusion $\overline{S(\Omega,n)} \supseteq \mathcal{S}$ is obvious, because the closure of a union contains the union of the closures. For the reverse inclusion, let $M \in \overline{S(\Omega,n)}$. Then, there exists a sequence of pencils $(M_i)_{i\in \N}$ such that $M_i \in S(\Omega,n)$ for all $i$ and $M = \lim_{i \rightarrow \infty} M_i$. For each $M_i$, there exists a generalized Schur form $M_i = Q_i T_i Z_i$, where $Q_i, Z_i \in U(n)$ and $T_i \in S_{tri}(\Omega,n)$. 
    
    For every $\epsilon > 0$, let $B_\epsilon$ denote the closed ball centered at $0$ and having radius $\|M\|_F+\epsilon$. Clearly, $K_\epsilon:=B_\epsilon \cap \overline{S_{tri}(\Omega,n)}$ is compact for all $\epsilon > 0$. As it holds that 
      
    \[\left|\| T_i \|_F - \| M \|_F\right| = \left|\| M_i \|_F - \| M \|_F\right| \leq \| M_i - M \|_F \quad \mathrm{and} \quad \lim_{i \rightarrow \infty} \| M_i - M \|_F = 0,\] there exists an index $j$ such that $(T_i)_{i>j}$ is a sequence in the compact set $K_\epsilon$. Moreover, the unitary group $U(n)$ is compact, and so is the Cartesian product $U(n) \times K_\epsilon \times U(n)$. As such, there exists a convergent subsequence $(Q_{f(i)},T_{f(i)},Z_{f(i)})_{i \in \N}$, where $f: \N \rightarrow \N$ is some strictly increasing function. Let us define $Q = \lim_{i \rightarrow \infty} Q_{f(i)},$ $T = \lim_{i \rightarrow \infty} T_{f(i)}$ and $Z = \lim_{i \rightarrow \infty} Z_{f(i)}$. It holds that
    \begin{align*}
        M = \lim_{i \rightarrow \infty} M_i = \lim_{i \rightarrow \infty} M_{f(i)} = \lim_{i \rightarrow \infty} Q_{f(i)} T_{f(i)} Z_{f(i)} = \lim_{i \rightarrow \infty} Q_{f(i)} \lim_{i \rightarrow \infty} T_{f(i)} \lim_{i \rightarrow \infty} Z_{f(i)} = Q T Z,
    \end{align*}
    where $Q, Z \in U(n)$ and $T \in K_\epsilon \subset \overline{S_{tri}(\Omega,n)}$. In other words, every limit point of the set $S(\Omega,n)$ belongs to $Q \overline{S_{tri}(\Omega,n)} Z$ for some fixed $Q, Z \in U(n)$, and hence it belongs to $\mathcal{S}$. This proves the statement.    
\end{proof}
%{\color{blue}VN: I added some new remarks below, and corrected elsewhere a statement that I had previously made but was not fully correct.}

Lemma \ref{lemma:closures}  allows us to justify a claim that we made in Section \ref{sec:problem}, that is, if a pencil  $P \in \overline{S(\Omega,n)}$, but $P \notin S(\Omega,n)$, then $P$ must be singular. We prove this fact formally in Proposition \ref{prop:closuresing} below.

\begin{proposition}\label{prop:closuresing}
    If $A + x B \in \overline{S(\Omega,n)}\setminus S(\Omega,n)$, then $A + x B$ is singular.
\end{proposition}
\begin{proof}
Suppose first that $n=1$. Consider a sequence $a_k + x b_k \subset S(\Omega,1)$ converging to $a+xb \in \overline{S(\Omega,1)} \setminus S(\Omega,1)$. Then, for all $k$, $(a_k,b_k) \neq (0,0)$ and $-a_k/b_k \in \Omega$. Suppose $(a,b) \neq (0,0)$, then $-a/b \in \overline{\C}$ is well defined; however, $\Omega$ is closed in $\overline{\C}$, and thus $-a/b \in \Omega$ implying $a+xb \in S(\Omega,1)$, which contradicts the assumption. Hence, necessarily $(a,b)=(0,0)$ and therefore $\overline{S(\Omega,1)} = S(\Omega,1) \cup \{0+x0\}$.

    If $n \geq 2$, let $S + x T = Q(A + x B)Z $, where $Q,Z \in U(n)$, denote a generalized Schur form of $A+xB \in \overline{S(\Omega,n)}\setminus S(\Omega,n)$. Then, Lemma \ref{lemma:closures} implies that $S + xT \in \overline{S_{tri}(\Omega,n)}$ but $S + x T \notin S_{tri}(\Omega,n)$. Since the diagonal entries of any pencil in $S_{tri}(\Omega,n)$ must belong to $S(\Omega,1)$, it follows from the first part of the proof that at least one diagonal element of $S+xT$ is equal to $0+x0$, and  therefore $\det(S+xT)=0$ and $S+xT$ is a singular pencil. In turn, this implies the claim, because $\det (A+xB) = \det(Q^*)\det(Z^*)\det (S+xT) = 0$ and thus $A+xB$ is singular.
\end{proof}  
Proposition \ref{prop:closuresing} shows that the limit points of $S(\Omega,n)$ are either pencils in $S(\Omega,n)$ (and hence regular) or singular pencils. In Example \ref{exnew} below we discuss whether the singular pencils in $\overline{S(\Omega,n)}\setminus S(\Omega,n)$ can have eigenvalues, and if so what they can be. Actually, it suffices in this sense to study $\overline{S_{tri}(\Omega,n)}\setminus S_{tri}(\Omega,n)$, in view of Lemma \ref{lemma:schur}, Proposition \ref{prop:closuresing}, and Definition \ref{def:evals}.
Note that the arguments in the proof of Proposition \ref{prop:closuresing} also show that the diagonal entries of a pencil in $\overline{S_{tri}(\Omega,n)}$ belong either to $S(\Omega,1)$ or to $\{0+x0\}$.
\begin{example}\label{exnew}
A singular pencil in $\overline{S_{tri}(\Omega,n)}$ may have no eigenvalues at all, and in fact this is the most typical situation. For example, suppose that $\Omega=\{ z \in \C : |z| \leq 1 \}$ and $n=2$, then
\[    \begin{bmatrix}
    0 & x-2 \\
    0 & 2x-1
\end{bmatrix}   \in \overline{S_{tri}(\Omega,2)} \setminus S_{tri}(\Omega,2)      \]
has no eigenvalues according to Definition \ref{def:evals}.

    Even if the nonzero diagonal entries of a singular pencil in $\overline{S_{tri}(\Omega,n)}$ must belong to $S(\Omega,1)$, not necessarily the eigenvalues of such a pencil (if any) must belong to $\Omega$. For a counterexample, take an arbitrary closed set $\Omega$, $n=3$, and consider the pencil
\[ A+xB=\begin{bmatrix}
    0 & \alpha x+\beta & 0\\
    0 & a+b x & \gamma x + \delta\\
    0 & 0 & 0
\end{bmatrix}   \in \overline{S_{tri}(\Omega,3)} \setminus S_{tri}(\Omega,3).\]
Here, $(0,0) \neq (a,b)$ is such that $-a/b \in \Omega$, whereas $(0,0)\neq (\alpha,\beta) \in \C^2$ and $(0,0) \neq (\gamma,  \delta) \in \C^2$ are arbitrary parameters. It is easy to see that the eigenvalues of $A+xB$ are $-\beta/\alpha$ and $-\delta/\gamma$. It is clear that we can choose the parameters $\alpha,\beta,\gamma,\delta$ so that these eigenvalues can be any elements of $\overline{\C}$; in particular, they do not necessarily need to belong to $\Omega$.

However, if a singular pencil in $\overline{S_{tri}(\Omega,n)}$ has precisely one zero diagonal element, then all its eigenvalues (if any) must belong to $\Omega$. This is because the eigenvalues of such a pencil (if any) must also be roots of one of its nonzero diagonal entries; see the proof of \cite[Theorem 5.3]{singpencil} for more details. For example, we can again take $\Omega=\{ z \in \C : |z| \leq 1 \}$ and $n=2$, and slightly modify the first example above to obtain
\[    \begin{bmatrix}
    0 & x-\frac12 \\
    0 & 2x-1
\end{bmatrix}   \in \overline{S_{tri}(\Omega,2)} \setminus S_{tri}(\Omega,2)        \]
which has the eigenvalue $\frac12 \in \Omega$.
\end{example}

% \textcolor{red}{I think we can phrase the result below by using lifts, too. Lifts with certain properties guarantee that the stationary points are the same (if I remember correctly).}
We are now in the position to state and prove Theorem \ref{thm:obj_function} which is the main result of this section.
\begin{theorem} \label{thm:obj_function}
    Let $A+x B \in \C[x]_1^{n \times n}$. Let $f(Q,Z)$ be the function on $U(n) \times U(n)$ defined by \eqref{eq:newobjfun}.
    Then:
    \begin{enumerate}
        \item The optimization problems \eqref{eq:riproblem} and \eqref{eq:minnewobjfun} have the same minimum value;
        \item The pair of unitary matrices $(Q_0,Z_0)$ is a global (resp. local) minimizer for \eqref{eq:minnewobjfun} if and only if the pencil $Q_0^* \mathcal{T}(Q_0 A Z_0 + x Q_0 B Z_0) Z_0^*$ is a global (resp. local) minimizer for \eqref{eq:riproblem}.
    \end{enumerate}
\end{theorem}
\begin{proof}
    %Let $S_{tri}(\Omega,n)$ be the set of $\Omega$-stable upper triangular pencils. 
    According to Lemma \ref{lemma:closures}, $\overline{S(\Omega,n)} = \bigcup_{Q,Z \in U(n)} Q \overline{S_{tri}(\Omega,n)} Z$. Hence, we obtain
    \[ \min_{S+x T \in \overline{S(\Omega,n)}} \|(A-S)+x(B-T)\|_F^2  = \min_{Q,Z \in U(n)} \min_{X+x Y \in \overline{S_{tri}(\Omega,n)}}  \|(A-Q^*XZ^*)+x(B-Q^*YZ^*) \|_F^2 \]
    \[ = \min_{Q,Z \in U(n)} \min_{X+x Y \in \overline{ S_{tri}(\Omega,n)}}  \|(QAZ-X)+x(QBZ-Y) \|_F^2 = \min_{Q,Z \in U(n)} \ \| QAZ+x QBZ - \mathcal{T}(QAZ+x QBZ) \|^2_F,\]
    using Proposition \ref{prop:nearest_tri} in the last step. This equation directly proves item 1; and, by applying Proposition \ref{prop:nearest_tri} again, item 2 also follows.
\end{proof}

Theorem \ref{thm:obj_function} shows that we can minimize the objective function $f(Q,Z)$ over $U(n) \times U(n)$ to compute the nearest $\Omega$-stable pencil. The optimization can in principle be addressed by any algorithm that is capable of doing optimization on Riemannian manifolds. For this, we need a closed-form expression for the projection $p_\Omega$. In the next section, we will derive such a formula for the cases of Hurwitz and Schur stability.

\section{Formulae for projection}\label{sec:projformula}

%Fix a closed set $\Omega \subseteq \overline{\mathbb{C}}$ and a square pencil $A+ x B$. %The nearest pencil with all eigenvalues in $\Omega$ can be found via Riemannian optimization if we can find a formula for the projection $p_{\Omega}: \C[x]_1 \rightarrow \C[x]_1^{\Omega}$. 
In this section, we derive a formula for the projection $p_{\Omega}: \C[x]_1 \rightarrow \overline{S(\Omega,1)}$ in two important cases: Hurwitz stability, where $\Omega$ is the left Riemann hemisphere, and Schur stability, where $\Omega$ is the bottom Riemann hemisphere. In the following, we identify the elements of $\C[x]_1$ with the elements of $\C^2$ via the map $a+x b \mapsto (a,b)$.

\subsection{Hurwitz stability} \label{sec:hurwitz}

Let $\Omega:=\{z \in \C \ : \ \Re(z) \leq 0\} \cup \{ \infty \}$, that is, $\Omega$ is the left hemisphere of the Riemann sphere. In this case, the set we want to project onto becomes $\overline{S(\Omega,1)} = \left\{a+b x: -\Re(\frac{a}{b}) \leq 0, \ a \in \C, \ b \in \C \backslash \{0\} \right\} \cup \{ a + 0x : a \in \C\}$. A direct calculation shows that this corresponds to projecting $a=\left(a_{1}+i a_{2}, a_{3}+i a_{4}\right) \in \mathbb{C}^{2}$ onto $\left\{x_{1} x_{3}+x_{2} x_{4} \geq 0\right\}$, where we have identified the element $(a_{1}+i a_{2}) + x(a_{3}+i a_{4}) \in \C[x]_1$ with $\left(a_{1}+i a_{2}, a_{3}+i a_{4}\right) \in \mathbb{C}^{2}$. Clearly, if $a_{1} a_{3}+a_{2} a_{4} \geq 0$, the projection is $p_{\Omega}(a)=a$, so let us assume $a_{1} a_{3}+a_{2} a_{4}<0$.

We can approach this constrained minimization problem with Lagrange's method. This yields a system of equations for unknowns $x_{i}$ and $\lambda$ such that
\begin{align} \label{eq:lagrange_system}
    \left[\begin{array}{cccc}
            1       & 0       & \lambda & 0       \\
            0       & 1       & 0       & \lambda \\
            \lambda & 0       & 1       & 0       \\
            0       & \lambda & 0       & 1
        \end{array}\right]\left[\begin{array}{l}
            x_{1} \\
            x_{2} \\
            x_{3} \\
            x_{4}
        \end{array}\right]=\left[\begin{array}{l}
            a_{1} \\
            a_{2} \\
            a_{3} \\
            a_{4}
        \end{array}\right], \quad x_{1} x_{3}+x_{2} x_{4}=0.
\end{align}
For notational simplicity, we set
\begin{align} \label{eq:alpha}
    \alpha:=\frac{\|a\|_{F}^{2}}{2\left(a_{1} a_{3}+a_{2} a_{4}\right)} \leq -1.
\end{align}
We note at this point that the case $\alpha = -1$ is degenerate and assume for now that $\alpha <-1$. Solving the polynomial equations \eqref{eq:lagrange_system} yields
\begin{align*}
    \left[\begin{array}{l}
            x_{1} \\
            x_{2} \\
            x_{3} \\
            x_{4}
        \end{array}\right]=\left(1-\lambda^{2}\right)^{-1}\left(\left[\begin{array}{l}
            a_{1} \\
            a_{2} \\
            a_{3} \\
            a_{4}
        \end{array}\right]-\lambda\left[\begin{array}{l}
            a_{3} \\
            a_{4} \\
            a_{1} \\
            a_{2}
        \end{array}\right]\right), \quad \lambda^{2}-2 \alpha \lambda+1=0 \Rightarrow \lambda=\alpha \pm \sqrt{\alpha^{2}-1} < 0.
\end{align*}
The choice with the positive sign in the formula for $\lambda$ yields the minimum (this can be seen either from the bordered Hessian or by substituting into the function). A direct but tedious calculation shows that the distance from the point $a$ to the set $\overline{S(\Omega,1)}$ is
\begin{align}\label{eq:dist}
    q(a) & =\left|\frac{\lambda}{\lambda^{2}-1}\right| \sqrt{\left(1+\lambda^{2}\right)\|a\|_{F}^{2}-4 \lambda\left(a_{1} a_{3}+a_{2} a_{4}\right)} \nonumber \\
         & =\sqrt{(a_{1} a_{3}+a_{2} a_{4})\lambda},
\end{align}
where $\lambda=\alpha+\sqrt{\alpha^{2}-1}$. In complex form, the projection in this case is
$$
    p_{\Omega}(a)=\frac{1}{1-\lambda^{2}}\left(a_{1}-\lambda a_{3}+i\left(a_{2}-\lambda a_{4}\right), a_{3}-\lambda a_{1}+i\left(a_{4}-\lambda a_{2}\right)\right) .
$$
The analysis above breaks down if $\alpha=-1$, in which case $\lambda=-1$. This is equivalent to the case $a_{1}=-a_{3}$ and $a_{2}=-a_{4}$. Substituting these values to \eqref{eq:lagrange_system} and solving the resulting system shows that the distance is still $q(a) = \sqrt{a_{1}^{2}+a_{2}^{2}}$ %and one (non-unique) projection is $p_{\Omega}(a)=\left(0+z\left(a_{3}+a_{4} i\right)\right)$. 
and all possible points at this distance are $(a_{1} + b_1 + i (a_{2} + b_2), b_{1} + i b_{2})$, where $b_{1} + i b_{2}$ is any complex number on a circle of radius $\left|a_{1} + i a_{2}\right|/2$ centered around the point $-(a_{1} + i a_{2})/2$.

\subsection{Schur stability} \label{sec:schur}

Let $\Omega:=\{z \in \C \ : \ \left|z\right| \leq 1\}$, that is, $\Omega$ is the unit disc of the complex plane. In this case, projecting an element of $\C[x]_1$ onto the set $\overline{S(\Omega,1)}$ corresponds to projecting $(a,b) \in \C^2$ onto the set $\{\left|a\right| \leq \left|b\right|\}$. Clearly, this projection is the identity map when $|a| \leq |b|$. When $|a| > |b|$ and $b \neq 0$, we get the intuitive result $$p_{\Omega}((a,b))=\left(\frac{\left|a\right|+\left|b\right|}{2\left|a\right|} a, \frac{\left|a\right|+\left|b\right|}{2\left|b\right|} b\right)$$ via Lagrange multipliers. When $|a| > |b|$ and $b=0$, the projection is non-unique. In this case, all the projections are of the form $p_{\Omega}((a,0))=\frac{1}{2} \left( a, \hat a \right)$, where $\hat a$ is any complex number such that $| \hat a | = | a |$. In all of these cases, the distance from $a+\lambda b$ to the set $S(\Omega,1)$ is
\begin{align}\label{eq:dist_schur}
    q((a,b)) & = \frac{\sqrt{2}(\left|a\right|-\left|b\right|)}{2}.
\end{align}

\section{Riemannian optimization}\label{sec:optimization}
Our aim is to use the trust-region method \cite{ABG07} of Manopt \cite{manopt} to minimize the objective function \eqref{eq:newobjfun} over the Riemannian manifold $U(n) \times U(n)$. In \cite{singpencil},  Dopico, Noferini and Nyman used a similar Riemannian optimization approach over the manifold $U(n) \times U(n)$. In \cite[Subsection 3.2]{singpencil}, one can find a detailed account on the technicalities related to optimization over this particular manifold, including expressions for the tangent space, projection to the tangent space as well as the computation of the Riemannian Hessian from the Euclidean Hessian. In this section, we present the main ingredients for what is needed, and refer the reader to \cite[Subsection 3.2]{singpencil} for more details.

\subsection{Gradient}
In order to numerically minimize the objective function \eqref{eq:newobjfun} over the Riemannian manifold $U(n) \times U(n)$, we need to calculate its Riemannian gradient. The Riemannian gradient can be computed as the tangent space projection of the Euclidean gradient \cite[Proposition 3.61]{rmanifolds}. In this subsection, we derive a formula for the Euclidean gradient of the objective function \eqref{eq:newobjfun}.

Written explicitly, the objective function is
\begin{align*}
    f(Q,Z) & = \| QAZ +  x QBZ - \mathcal{T}(QAZ +  x QBZ) \|^2                                       \\
           & = \| QAZ - \mathcal{T}_0(QAZ +  x QBZ) \|^2 + \| QBZ - \mathcal{T}_1(QAZ +  x QBZ) \|^2.
\end{align*}
Both here and further below, we introduce the notation $\mathcal{T}_0$ and $\mathcal{T}_1$ to denote, respectively, the trailing and leading coefficients of $\mathcal{T}(A+xB)$, where $\mathcal{T}$ is defined as in Proposition \ref{prop:nearest_tri}. Moreover, we will use a similar notation, with subscripts $0$ and $1$, to denote the matrix coefficients of any pencil-valued function.

Let us first consider the gradient of $f$ with respect to the first argument $Q$. To this end, let $f = f_1 \circ f_2$, where $f_1(A + xB) = \| (A +  x B) - \mathcal{T}(A +  x B) \|^2$ and $f_2(Q) = QAZ +  x QBZ$. We use the following lemma, which is based on \cite[Theorem 2.1]{np}, which in turn generalizes the result in \cite[Proposition 4.1]{Chazal} to unbounded closed sets.
% \textcolor{red}{The above looks weird now when the references are in the superscript. Maybe we have to change the phrasing throughout the paper to account for the change in how the references are displayed.}
\begin{lemma} \label{lemma:gradient}
    The function $f_1(A + xB) = \| (A +  x B) - \mathcal{T}(A +  x B) \|^2$ is almost everywhere differentiable. Moreover, the points of non-differentiability coincide with the medial axis. Outside the medial axis, the gradient is given by $$\nabla f_1(A + xB) = 2(A +  x B - \mathcal{T}(A +  x B)).$$
\end{lemma}
\begin{proof}
    See \cite[Proposition 4.1]{Chazal} and \cite[Theorem 2.1]{np}.
\end{proof}
For the Hurwitz case, that is, the case $\Omega=\{z \in \C \ : \ \Re(z) \leq 0\} \cup \{ \infty \}$, the pencil $A+xB$ lies in the medial axis if and only if, for any of its diagonal entries that belong to $S(\Omega^c,1)$, it holds that $\alpha = -1$, where $\alpha$ is defined as in \eqref{eq:alpha}. For the Schur case, that is, the case $\Omega=\{z \in \C \ : \ \left|z\right| \leq 1\}$, the pencil $A+xB$ lies in the medial axis if and only if, for any of its diagonal elements $a_{ii} + x b_{ii}$, it holds that $a_{ii} \neq 0$ and $b_{ii}=0$. By Lemma \ref{lemma:gradient}, the distance function $f_1$ is differentiable everywhere except at these points. 

Outside the medial axis, we can use Lemma \ref{lemma:gradient} to express the directional derivative of the function $f_1$ at the point $A+x B$ in the direction $E + xF$ as
\begin{align*}
    D f_1(A+x B)[E + x F] = 2 \langle (A +  x B) - \mathcal{T}(A +  x B), E + x F\rangle,
\end{align*}
where $\langle \cdot , \cdot \rangle$ is the real inner product in $\C^{n \times n} \times \C^{n \times n}$ defined as
\begin{equation*} \label{eq:innerprodunun}
    \langle A_1 + x A_2 \, , \, B_1 + x B_2 \rangle := \mathrm{Re} \left(\mbox{trace} (A_1^* B_1) + \mbox{trace} (A_2^* B_2) \right),
\end{equation*}
for $ (A_1, A_2) ,  (B_1, B_2) \in \C^{n \times n} \times \C^{n \times n}$. Similarly, for two constant matrices $A,B \in \C^{n \times n}$ we define $\langle A, B \rangle := \mathrm{Re} \left(\mbox{trace} (A^* B)  \right)$. Moreover, clearly the directional derivative of the function $f_2$ at the point $Q$ in the direction $E$ is
\begin{align*}
    D f_2(Q)[E] = EAZ + x EBZ.
\end{align*}
Thus, we can express the directional derivative of the composition $f = f_1 \circ f_2$ at the point $Q$ in the direction $E$ as
\begin{align*}
    D f(Q)[E] & = (D f_1(f_2(Q)) \circ D f_2(Q) ) [E]                                       \\
              & = 2 \langle (QAZ +  x QBZ) - \mathcal{T}(QAZ +  x QBZ), EAZ + x EBZ \rangle \\
    % % &= 2 \langle (QAZ +  x QBZ) - \mathcal{T}(QAZ +  x QBZ), EAZ\rangle \\
    % % & \ + 2 \langle (QAZ +  x QBZ) - \mathcal{T}(QAZ +  x QBZ), QAF\rangle \\
    % % & \ + 2 \langle (QAZ +  x QBZ) - \mathcal{T}(QAZ +  x QBZ), x EBZ\rangle \\
    % % & \ + 2 \langle (QAZ +  x QBZ) - \mathcal{T}(QAZ +  x QBZ), x QBF\rangle \\
    %           & = 2 \tr ((QAZ - \mathcal{T}_0 (QAZ +  x QBZ))^* EAZ)                        \\
    % % & \ + 2 \tr ((QAZ - \mathcal{T}_0 (QAZ +  x QBZ))^* QAF) \\
    %           & \ + 2 \tr ((QBZ - \mathcal{T}_1(QAZ +  x QBZ))^* EBZ)                       \\
    % % & \ + 2 \tr ((QBZ - \mathcal{T}_1(QAZ +  x QBZ))^* QBF) \\
    %           & = 2 \tr (((QAZ - \mathcal{T}_0 (QAZ +  x QBZ))(AZ)^*)^* E)                  \\
    % % & \ + 2 \tr (((QA)^*(QAZ - \mathcal{T}_0 (QAZ +  x QBZ)))^* F) \\
    %           & \ + 2 \tr (((QBZ - \mathcal{T}_1(QAZ +  x QBZ))(BZ)^*)^* E)                 \\
    % % & \ + 2 \tr (((QB)^*(QBZ - \mathcal{T}_1(QAZ +  x QBZ)))^* F) 
              & = \langle 2 (QAZ - \mathcal{T}_0 (QAZ +  x QBZ))(AZ)^*, E \rangle           \\
              & \ + \langle 2 (QBZ - \mathcal{T}_1(QAZ +  x QBZ))(BZ)^*, E \rangle.
\end{align*}
The gradient is the unique vector $v$ with the property that $\langle v, w \rangle$ yields the directional derivative in the direction $w$ for all $w$. This leads us to the following formula for the gradient:
\begin{align} \label{eq:gradient}
    \nabla f_Q(Q) = 2 (QAZ - \mathcal{T}_0 (QAZ +  x QBZ))(AZ)^* + 2 x (QBZ - \mathcal{T}_1(QAZ +  x QBZ))(BZ)^*.
\end{align}
The gradient with respect to the second argument $Z$ can be computed analogously.

Note that there is no need to compute the gradient of the operator $\mathcal{T}$ or the projection $p_{\Omega}$ hidden inside it. This is a consequence of Lemma \ref{lemma:gradient}. In this way, the formula $\eqref{eq:gradient}$ is fully general and works for any underlying closed set $\Omega$.

\subsection{Hessian}
% We partition the gradient with respect to $Q$ into two parts as follows:
% \begin{align*}
%     &\nabla f_Q(Q) = 2 L(QAZ)(AZ)^* + 2 L(QBZ)(BZ)^* \\
%     &+ 2 (\text{Diag}(QAZ - \mathcal{T}_0 (QAZ +  x QBZ)))(AZ)^* + 2 (\text{Diag}(QBZ - \mathcal{T}_1(QAZ +  x QBZ)))(BZ)^*.
% \end{align*}
% where $L$ is the operator that sets the upper triangular part of the matrix to be zero, and $\text{Diag}$ is the operator that sets every element to zero except the diagonal elements.

% The directional derivative of the first part with respect to $Q$ is 
% \begin{align*}
%     \mathrm{D} \nabla_{Q} g (Q,Z)[d_Q,d_Z] = 2 L(d_Q AZ)(AZ)^* + 2 L(QA d_Z)(AZ)^* + 2 L(QAZ)(A d_Z)^*.
% \end{align*}

In order to use a second-order method for the minimization of \eqref{eq:newobjfun}, in this subsection we derive a formula for the Euclidean Hessian of the objective function. We refer the reader to \cite[Subsection 3.2]{singpencil} for more details on how to construct the Riemannian Hessian from the Euclidean Hessian.
% In order to use a second-order method for the minimization of \eqref{eq:newobjfun}, in this section we derive a formula for the Euclidean Hessian of the objective function. 
The Euclidean Hessian is defined via the directional derivative of the gradient \cite[p. 23]{rmanifolds}. 

We start by noting that the objective function $\eqref{eq:newobjfun}$ is twice differentiable almost everywhere. In particular, the points where the function might not be twice differentiable either lie on the medial axis or belong to $\partial{S(\Omega,n)}$, that is, the boundary of the set $S(\Omega,n)$. In fact, the function is smooth outside this set of zero measure. In this subsection, we tacitly assume that the Hessian is only evaluated at points where the function is twice differentiable.

Now, let $L: \C^{n\times n} \rightarrow \C^{n\times n}$ denote the linear operator that sets the strictly upper triangular part of its input matrix to be zero. Then, the directional derivative of the gradient  $\eqref{eq:gradient}$ with respect to the first argument $Q$ in the direction $(d_Q,d_Z) \in \C^{n \times n} \times \C^{n \times n}$ is
\begin{equation}
    \begin{aligned} \label{eq:hessian}
         & D \nabla f_Q(Q)[d_Q,d_Z] =                       \\
         & 2 (QAZ - \mathcal{T}_0 (QAZ +  x QBZ))(A d_Z)^*  \\
         & + 2 (L(d_Q AZ + QA d_Z) - H_0)(AZ)^*             \\
         & + 2 (QBZ - \mathcal{T}_1(QAZ +  x QBZ))(B d_Z)^* \\
         & + 2 (L(d_Q BZ + QB d_Z) - H_1)(BZ)^*,
    \end{aligned}
\end{equation}
where $d_Q$ and $d_Z$ are the matrix directions for the matrices $Q$ and $Z$, respectively, and $H_0$ and $H_1$ are diagonal matrices s.t.
\begin{align*}
     & (H_0)_{ii} = (D h_{i}(Q,Z)[d_Q,d_Z])_0, \\
     & (H_1)_{ii} = (D h_{i}(Q,Z)[d_Q,d_Z])_1.
\end{align*}
Here, $h_i(Q,Z) := p_{\Omega}((QAZ + xQBZ)_{ii})$. A similar result holds for the gradient with respect to the second argument $Z$.

In order to compute \eqref{eq:hessian} (in particular the matrices $H_0$ and $H_1$), we need to compute the directional derivatives of the functions $h_i$, which involves computing the directional derivative of the projection $p_{\Omega}$. We do this for the projections associated with Hurwitz stability and Schur stability.

\subsubsection{Hurwitz stability}
As seen in Subsection \ref{sec:hurwitz}, the projection in the region $a_{1} a_{3}+a_{2} a_{4}<0$ and outside the medial axis for Hurwitz stability is
$$
    p_{\Omega}(a)=\frac{1}{1-\lambda^{2}}\left(a_{1}-\lambda a_{3}+i\left(a_{2}-\lambda a_{4}\right), a_{3}-\lambda a_{1}+i\left(a_{4}-\lambda a_{2}\right)\right) .
$$
In the following, we identify elements in $\C[x]_1$ and $\mathbb{C}^{2}$ with elements in $\R^4$ via the maps $z_1 + i z_2 + x (z_3 + i z_4) \mapsto (z_1,z_2,z_3,z_4)$ and $\left(z_{1}+i z_{2}, z_{3}+i z_{4}\right) \mapsto (z_1,z_2,z_3,z_4)$, respectively. 

A direct calculation shows that 
\begin{align*}
    \frac{d \lambda}{d \alpha} & = 1 + \frac{\alpha}{\sqrt{\alpha^2 -1}},                                                                                                  \\
    D \lambda                  & = \frac{d \lambda}{d \alpha} D \alpha,                                                                                                    \\
    D \alpha(a)[b]             & = \frac{a^T b}{\left(a_{1} a_{3}+a_{2} a_{4}\right)} - \frac{\left|a\right|_{F}^{2}}{2\left(a_{1} a_{3}+a_{2} a_{4}\right)^2} [a_3,a_4,a_1,a_2] b,
\end{align*}
where $\lambda$ and $\alpha$ are defined as in Subsection \ref{sec:hurwitz}. When $a_{1} a_{3}+a_{2} a_{4}<0$, the directional derivative of the projection is
\begin{align*}
    D p_{\Omega}(a)[b]= & \frac{2 \lambda}{(1-\lambda^{2})^2} D \lambda(a)[b] (b_{1}-\lambda b_{3}+i\left(b_{2}-\lambda b_{4}\right) - D \lambda(a)[b](a_3 + i a_4), \\
               & b_{3}-\lambda b_{1}+i\left(b_{4}-\lambda b_{2}\right) - D \lambda(a)[b](a_1 + i a_2))
    % =&  \frac{2 \lambda}{1-\lambda^{2}} D \lambda(a)[b] P(b) - \frac{2 \lambda}{(1-\lambda^{2})^2} (D \lambda(a)[b])^2 (a_3 + i a_4, a_1 + i a_2) \\
    % =&  \frac{2 \lambda}{1-\lambda^{2}} D \lambda(a)[b] (P(b) - \frac{D \lambda(a)[b]}{1-\lambda^{2}}) (a_3 + i a_4, a_1 + i a_2)
    .
\end{align*}
Let us now fix the index $i$ and define the following functions: 
\begin{align*}
&g_i(Q,Z) := (QAZ + xQBZ)_{ii},\\
&h_i(Q,Z) := (p_{\Omega} \circ g_i)(Q,Z) = p_{\Omega}((QAZ + xQBZ)_{ii}). 
\end{align*}
Note that $h_i$ is the function whose directional derivative we aim to compute, as it is required to compute the Hessian \eqref{eq:hessian}. We begin by observing that 
\begin{align*}
    D g_i(Q,Z)[d_Q,d_Z] = (d_Q AZ + xd_Q BZ + Q A d_Z + x Q B d_Z)_{ii}.
\end{align*}
Then, the directional derivative of $h_i$ in the direction $(d_Q,d_Z)$ is
\begin{align*}
     & D h_{i}(Q,Z)[d_Q,d_Z] = (D p_{\Omega}(g_i(Q,Z)) \circ D g_i(Q,Z)) [d_Q,d_Z]                                                                            \\
     & = \frac{2 \lambda}{(1-\lambda^{2})^2} D \lambda(g_i(Q,Z))[D g_i(Q,Z)[d_Q,d_Z]] ( (QAZ - \lambda QBZ)_{ii}, (QBZ - \lambda QAZ)_{ii})          \\
     & + \frac{1}{1-\lambda^{2}} ( (d_Q AZ - \lambda d_Q BZ + Q A d_Z - \lambda Q B d_Z)_{ii} - D \lambda(g_i(Q,Z))[D g_i(Q,Z)[d_Q,d_Z]] (QBZ)_{ii}, \\
     & (d_Q BZ - \lambda d_Q AZ + Q B d_Z - \lambda Q A d_Z)_{ii} - D \lambda(g_i(Q,Z))[D g_i(Q,Z)[d_Q,d_Z]] (QAZ)_{ii}),
\end{align*}
where
\begin{align*}
     & D \lambda(g_i(Q,Z))[D g_i(Q,Z)[d_Q,d_Z]] = (1 + \frac{\alpha}{\sqrt{\alpha^2 -1}})                                                              \\
     & \cdot (\frac{s^T t}{\left(s_{1} s_{3}+s_{2} s_{4}\right)} - \frac{\|s\|_{F}^{2}}{2\left(s_{1} s_{3}+s_{2} s_{4}\right)^2} [s_3,s_4,s_1,s_2] t),
\end{align*}
and
\begin{align*}
     & s = (QAZ + x QBZ )_{ii},                             \\
     & t = (d_Q AZ + xd_Q BZ + Q A d_Z + x Q B d_Z)_{ii}.
\end{align*}
The above holds when $s_{1} s_{3}+s_{2} s_{4}<0$, otherwise $s \in \overline{S(\Omega,1)}$ and the directional derivative of $h_i$ is zero.

% A similar result holds for the directional derivative of the gradient with respect to the second argument $Z$.

\subsubsection{Schur stability}
% The same formula as in the previous section holds if the matrices $H_0$ and $H_1$ are modified appropriately. These matrices are again defined similarly via the directional derivative of $h_i(Q,Z) := (P \circ g_i)(Q,Z) = P((QAZ + xQBZ)_{ii})$, where

As seen in Subsection \ref{sec:schur}, the projection for the Schur case is $p_{\Omega}((a,b))=\left(\frac{\left|a\right|+\left|b\right|}{2\left|a\right|} a, \frac{\left|a\right|+\left|b\right|}{2\left|b\right|} b \right)$ outside the medial axis and in the region $|a| > |b|$ and $b \neq 0$. In this case, the directional derivative becomes
\begin{align*}
    D p_{\Omega}((a,b))[(c,d)]= & \left(\frac{1}{2}\left(c + \frac{\langle b, d \rangle}{\left|a\right| \left|b\right|} a - \frac{\langle a, c \rangle \left|b\right|}{\left|a\right|^3} a + \frac{\left|b\right|}{\left|a\right|} c\right), \right. \\
                       & \left. \frac{1}{2}\left(d + \frac{\langle a, c \rangle}{\left|a\right| \left|b\right|} b - \frac{\langle b, d \rangle \left|a\right|}{\left|b\right|^3} b + \frac{\left|a\right|}{\left|b\right|} d\right)\right),
\end{align*}
where $\langle \cdot, \cdot \rangle$ is the real inner product defined as $\langle a, b \rangle := \text{Re}(a^* b)$. Let us now fix the index $i$ and consider the functions
\begin{align*}
&g_i(Q,Z) := ((QAZ)_{ii},(QBZ)_{ii}),\\
&h_i(Q,Z) := (p_{\Omega} \circ g_i)(Q,Z) = p_{\Omega}(((QAZ)_{ii},(QBZ)_{ii})). 
\end{align*}
Then,
\begin{align*}
     & D h_{i}(Q,Z)[d_Q,d_Z] = (D p_{\Omega}(g_i(Q,Z)) \circ D g_i(Q,Z)) [d_Q,d_Z]                                                                                                 \\
     & = \left(\frac{1}{2}\left(t_1 + \frac{\langle s_2, t_2 \rangle}{\left|s_1\right| \left|s_2\right|} s_1 - \frac{\langle s_1, t_1 \rangle \left|s_2\right|}{\left|s_1\right|^3} s_1 + \frac{\left|s_2\right|}{\left|s_1\right|} t_1\right), \right. \\
     & \left. \frac{1}{2}\left(t_2 + \frac{\langle s_1, t_1 \rangle}{\left|s_1\right| \left|s_2\right|} s_2 - \frac{\langle s_2, t_2 \rangle \left|s_1\right|}{\left|s_2\right|^3} s_2 + \frac{\left|s_1\right|}{\left|s_2\right|} t_2\right)\right),
\end{align*}
where
\begin{align*}
     & s_1 = (QAZ)_{ii},              \\
     & s_2 = (QBZ)_{ii},              \\
     & t_1 = (d_Q AZ + Q A d_Z)_{ii},  \\
     & t_2 = (d_Q BZ + Q B d_Z)_{ii}.
\end{align*}
The above holds when $|s_1| > |s_2|$; otherwise $(s_1 + x s_2) \in \overline{S(\Omega,1)}$ and the directional derivative of $h_i$ is zero.

\section{Numerical experiments}\label{sec:numexp}
In this section, we perform numerical experiments and benchmark the performance of our algorithm against other existing methods. To our knowledge, the method presented by Gillis, Mehrmann and Sharma \cite{gms} is the only existing algorithm devised for the problem of computing the distance to stability of a pencil. 

The code for the method presented in this paper can be found in the GitHub repository \href{https://github.com/NymanLauri/nearest-stable-pencil}{https://github.com/NymanLauri/nearest-stable-pencil}. It requires the Manopt package \cite{manopt} for optimization over Riemannian manifolds, which can be downloaded from \href{https://www.manopt.org/}{manopt.org}. Numerical tests were performed by using MATLAB R2023a and a modified version of Manopt 7.1. The one slight modification, with the goal to remove an unnecessary bottleneck, was that the check on line 36 in file \texttt{multihconj.m} was removed. This line checks whether the file \texttt{pagectranspose} exists, which we know exists for MATLAB R2020b and newer. 

This section is organized as follows. First, in Subsection \ref{sec:numhurwitz}, we perform numerical experiments for the case of Hurwitz stability. In this subsection, we provide heuristic observations and compare our method against other existing methods in the literature. More specifically, we compare against the DH method outlined in the work of Gillis, Mehrmann and Sharma\cite{gms}, which is a specialized algorithm for finding the nearest real Hurwitz stable pencil. We both compare against the DH method for individual examples and present a comparison based on statistical experiments. Then, in Subsection \ref{sec:numschur}, we test the implementation for the case of Schur stability. We do this by running the algorithm for the same individual examples as with the Hurwitz case.  

% we provide heuristic observations about the output of our algorithm. Then, in subsections \ref{sec:individual} and \ref{sec:statistical} we compare our method against other existing methods in the literature. More specifically, we compare against the DH method outlined in the work of Gillis, Mehrmann and Sharma\cite{gms}, which is a specialized algorithm for finding the nearest real Hurwitz stable pencil. In subsection \ref{sec:individual}, we compare against the DH method for individual examples, and afterwards in subsection \ref{sec:statistical}, we present a comparison based on statistical experiments. Finally, in subsection \ref{sec:schur_experiments} we run some experiments for the Schur case.

% We list some heuristic observations and benchmark the method against  perform that of \cite{gms}. 

\subsection{Hurwitz stability} \label{sec:numhurwitz}

% In this subsection, we list some heuristic observations about the output of our algorithm, in the case of Hurwitz stability. More rigorous statistical experiments are presented in subsection \ref{sec:statistical}.

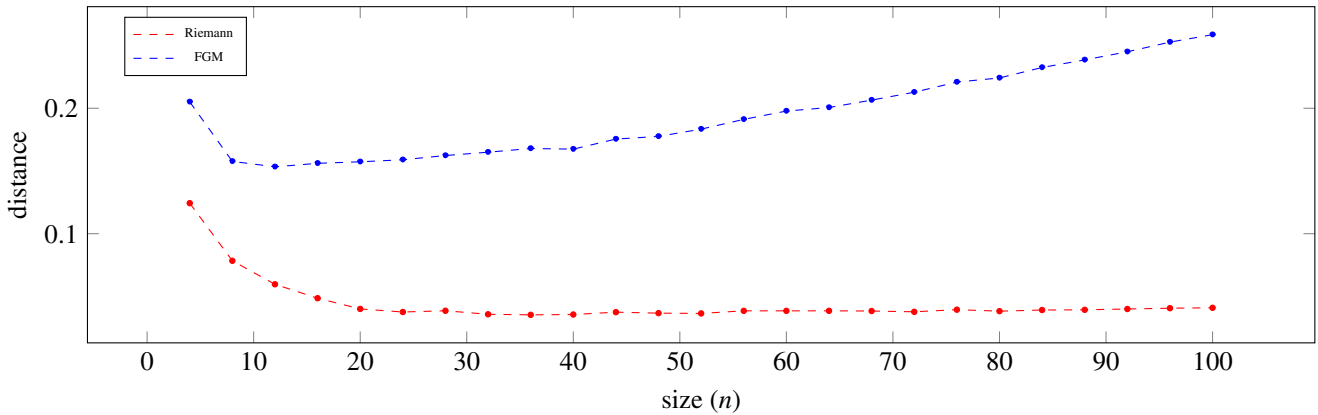
\begin{figure}[h!]
    \centering
    \begin{tikzpicture}
        \begin{axis}[
                legend pos = north west,
                width=\linewidth,
               height= 6cm,
                xlabel={size ($n$)},
                ylabel={distance}
            ]
             \addplot[red, dashed] table[y index = 1] 
            {statistical_comp.dat};
             \addplot[blue, dashed] table[y index = 2] 
            {statistical_comp.dat};
            \addplot[only marks, red, mark size =1pt] table[y index = 1] 
            {statistical_comp.dat};
            \addplot[only marks, blue, dashed, mark size =1pt] table[y index = 2] 
              {statistical_comp.dat};
            \legend{{\small{Riemann}}, {\small{FGM}}}
          \end{axis}
        \end{tikzpicture}

    % \begin{tikzpicture}
    % \begin{axis}[
    %         legend pos = north east,
    %         width=\linewidth,
    %         % xtick={5,10, 15, 20, 25, 30},
    %         %xticklabels={},
    %       % ymax=100,
    %        height= 6cm,
    %        title={{\small{Median running time}}},
    %         xlabel={size ($n$)},
    %         ylabel={time (s)}
    %     ]
    %      \addplot[red, dashed] table[y index = 3] 
    %     {statistical_comp.dat};
    %      \addplot[blue, dashed] table[y index = 4] 
    %     {statistical_comp.dat};
    %     \addplot[only marks, red, mark size =1pt] table[y index = 3] 
    %     {statistical_comp.dat};
    %     \addplot[only marks, blue, dashed, mark size =1pt] table[y index = 4] 
    %       {statistical_comp.dat};
    %     \legend{{\small{Riemann}}, {\small{DH}}}
    %   \end{axis}
    % \end{tikzpicture}

    \caption{Statistical comparison between the approach of this paper (labelled as Riemann in the picture) and the FGM algorithm \cite{gms} for $4 \leq n \leq 100$.}
    \label{fig:DH_comparison}
\end{figure}

We start by providing some general remarks. 
 Numerical experiments with randomly generated pencils (entries drawn from the normal distribution) suggest that, typically:
\begin{enumerate}
    \item The computed pencil has a long Jordan chain. For example, a pencil of size 10 might have Jordan chains of length (7,2,1), or something similar. To make this statement more precise, we ran the algorithm for a sample of 100 input of sizes 5 and 10. For size 5, 80 of the outputs had a nontrivial Jordan chain, while for size 10, 89 of the outputs had a nontrivial Jordan chain. In this experiment, the algebraic multiplicity of an eigenvalue was calculated using a tolerance of $10^{-12}$, while the numerical rank of a matrix was determined by the number of singular values below the tolerance $10^{-12}$. If the algebraic multiplicity of an eigenvalue $\lambda$ was higher that the numerical rank of the associated matrix $A + \lambda B$, the pencil $A+xB$ was considered to have a nontrivial Jordan chain.
    \item The clear majority of the eigenvalues have real part equal to zero. %Sometimes all eigenvalues have real part equal to zero. 
    Sometimes many of the eigenvalues are zero or infinity.
    \item For an input of size 10, an $\Omega$-stable pencil is found at a distance of approximately 5\% of the norm of the input, which seems reasonable (for more details, see Figure \ref{fig:DH_comparison}).
    \item The above observations seem to hold for both real and complex input.
\end{enumerate}
Given that the input are Gaussian, the first and second observations may at first sound counterintuitive and perhaps even suspicious. However, these results are analogous to those in \cite{np}, where such minimizers were also observed frequently in numerical experiments. In \cite{np}, it was also rigorously shown that, for a randomly generated (with an absolutely continuous distribution) input, the nearest stable matrix has a nontrivial Jordan chain with nonzero probability; see \cite[Example 5.6]{np}. We conjecture that the result can be extended to the pencil setting.

Next, we compare our method against the DH method outlined in the work by Gillis, Mehrmann and Sharma \cite{gms}, which is a specialized algorithm for finding the nearest real Hurwitz stable pencil. Even though our method and the DH method solve different variants of the problem (see below), we can still perform comparisons between these methods. In particular, it is possible to naively use our method on the real case by simply optimizing over the orthogonal group as opposed to the unitary group. This guarantees that the output will be real; on the other hand, doing it in this manner entirely excludes potential minimizers that have complex eigenvalues (which is a set of non-zero measure). This is potentially a disadvantage for our method, and it may be overcome by designing a real variant of the algorithm, which is a possible topic of future research but beyond the scope of this paper. However, for the sake of comparison and keeping in mind that this choice potentially disadvantages our approach but not the DH method, in the experiment we adopt this naive approach.

The DH method attempts to find the nearest pencil that is regular, of index at most one (i.e. no Jordan chains longer than one exist for infinite eigenvalues), and has all finite eigenvalues in the open left half plane. The main difference with our approach is that we allow for arbitrarily long Jordan chains for infinite eigenvalues. Recall the remark above that, according to our statistical experiments on randomly generated input, in the minimizer long Jordan chains for both finite and infinite eigenvalues can occur in practice, and indeed seem to appear rather frequently.

Before presenting a comparison, we note in Lemma \ref{lem:samedistance} below that the peculiarities in the feasible set for the DH method cannot actually affect the distance to be computed, and therefore the comparison is fair. For simplicity of exposition, in Lemma \ref{lem:samedistance} we do not consider the real constraint, but even in that case the situation is the same with a similar proof (except that a real generalized Schur form must be used).

\begin{lemma}\label{lem:samedistance}
Let  $\Omega \subset \overline{\C}$ be the closed left half hemisphere of the Riemann sphere, and let $\mathcal{H} \subset \C[x]_1^{n \times n}$ denote the set of regular pencils of index at most $1$ and whose eigenvalues are all in $\Omega$. Then, $\overline{\mathcal{H}} = \overline{S(\Omega,n)}$. In particular, for every pencil $X \in \C[x]_1^{n \times n}$, it holds $d(X,\mathcal{H})=d(X,S(\Omega,n)).$
    \end{lemma}
    \begin{proof}
Since $\mathcal{H} \subseteq S(\Omega,n)$, it is obvious that $\overline{\mathcal{H}} \subseteq \overline{S(\Omega,n)}$. Therefore, we just need to prove $ \overline{S(\Omega,n)} \subseteq \overline{\mathcal{H}}$.

For $n=1$, $\mathcal{H} = S(\Omega,n)$ and the result follows trivially. For $n\geq2$, note that $\mathcal{H} \subsetneq S(\Omega,n)$, and denote by $\mathcal{J} $ the set of pencils that belong to $S(\Omega,n)$ but not to $\mathcal{H}$. Then, $\mathcal{J}$ contains precisely the pencils that are regular, with eigenvalues all in $\Omega$, and with index at least $2$. It is clear that the statement is implied by the following claim: For every pencil $P \in \mathcal{J}$ and every $\varepsilon > 0$, there exists $P' \in \mathcal{H}$ with $\|P-P'\|_F=\epsilon$. Indeed, if $Y \notin \mathcal{H}$ is a limit point of $S(\Omega,n)$, then for all $\epsilon > 0$ there is $P \in S(\Omega,n)$ such that $0 <\|Y-P\|\leq \epsilon$, and hence (by taking $P'=P$ if $P \in \mathcal{H}$ and $P'$ as above if $P \in \mathcal{J}$) there exists $P' \in \mathcal{H}$ with $0<\|Y-P'\| \leq 2 \epsilon$, thus implying that $Y$ is also a limit point of $\mathcal{H}$.

        It remains to prove the claim. To this goal, let $P=QTZ$ be a generalized Schur form of $P$, and suppose that $P$ has $k \geq 2$ infinite eigenvalues. Then there are precisely $k$ diagonal elements in the pencil $T$ of the form $\rho_i e^{i \alpha_i}+x0$, where, for all $i=1,\dots,k$, $\rho_i>0$ and $\alpha_i \in [-\pi,\pi]$. Let $\epsilon>0$. For $i=2,\dots,k$ we can perturb the $i$-th such element of $T$ to $e^{i \alpha_i} (\rho_i + x \epsilon/\sqrt{k-1})$, thus moving the eigenvalue to $-\sqrt{k-1}\rho_i/\epsilon \in \Omega$. Thus, the resulting perturbed pencil $P'$ has a simple infinite eigenvalue, while all the other eigenvalues are finite and in $\Omega$; therefore, $P' \in \mathcal{H}$. Moreover, such perturbation has norm equal to $\|P-P'\|_F=\epsilon$.
    \end{proof}

% Equivalently, we can start from a real point (e.g., the identity), and we'll stay real. 

% \textcolor{red}{For real pencils, the image of a gradient flow that starts from a real point is always a subset of $O(n)^2$. This would lead to a real flow and a real minimizer. This is natural behaviour of our algorithm: we can simply say that our algorithm "finds a real minimizer", and we can compare this minimum to the minimum found by \cite{gms}. Of course, it is still evident that in this way, we cannot get complex eigenvalues.}

We perform the comparison by running our algorithm on the same example problems that were originally considered by Gillis, Mehrmann and Sharma in\cite{gms}. For the comparisons, we used the implementation ``Nearest Stable Matrix Pairs" v1.0 that can, at the time of writing, be accessed via the webpage of Nicolas Gillis: \href{https://sites.google.com/site/nicolasgillis/code}{https://sites.google.com/site/nicolasgillis/code}.

\begin{example} \label{ex1}
    Let $n = 20$ and consider the pencil $I_n x - M$, where M is a banded Toeplitz matrix, called the \texttt{grcar} matrix, in which the main diagonal and three superdiagonals are set to 1, while the first subdiagonal is set to -1.

    Gillis, Mehrmann and Sharma\cite{gms} report finding a real Hurwitz stable pencil at the squared distance $6.28$. Using a randomly generated point as a starting point, our algorithm finds a real Hurwitz stable pencil at the squared distance $1.99$. This pencil is saved in the \texttt{example6dot1.mat} file that can be found in the GitHub repository.\footnote{The coefficients of the computed pencil $Tx + S$ are saved in the variables $S$ and $T$. The starting point is saved in the $x0$ variable, and the variable $Q$ stores the orthogonal transformations that triangularize the pencil. We note that computing the eigenvalues in this basis is much more accurate than calling \texttt{eig(S,-T)} directly.} It is worth noting that the computed pencil has a Jordan chain of length five corresponding to an infinite eigenvalue, and hence it is not included in the feasible set of the DH method; nevertheless, by Lemma \ref{lem:samedistance} it belongs to the closure of such feasible set. Therefore, the squared distance to the feasible set of the DH method cannot be any higher than $\approx 1.99$.
    
    %It is worth noting that this pencil has a finite eigenvalue with multiplicity 20.
    %     \[
    % M = \begin{bmatrix}
    %      1 &  1 & 1  & 1  &   &   &  & \\
    %     -1 &  1 & 1 &  1 &  1 &  &  & \\
    %       & -1 & 1 & 1 & 1 & 1 &  & \\
    %       &  & -1 & 1 & 1 & 1 & 1 & \\
    %       &  &  & -1 & 1 & 1 & 1 & 1\\
    %       &  &  &  & -1 & 1 & 1 & 1\\
    %       &  &  &  &  & -1 & 1 & 1\\
    %       &  &  &  &  &  & -1 & 1
    % \end{bmatrix},
    % \]
\end{example}

\begin{example} \label{ex2}
    Consider a damped oscillating system which can be modelled with a differential equation of form $M \ddot x + C \dot x + K x = f(t) $. We set $m = c = k = [ 1 \ 2 \ 3 \ 4 \ 5 \ 6 \ 7 \ 8 \ 9 \ 10 ] $ (see \cite{veselic} for how the matrices $M, C, K$ are constructed).

    Let $n=10$ and $\epsilon = 0.1$. We consider the pencil $A + x B$, where
    $$
        B=\left[\begin{array}{cc}
                M & 0   \\
                0 & I_n
            \end{array}\right], A=-(J-R) Q, \quad J=\left[\begin{array}{cc}
                0   & -I_n \\
                I_n & 0
            \end{array}\right], \quad R=\left[\begin{array}{cc}
                C & 0             \\
                0 & -\epsilon I_n
            \end{array}\right], \quad Q=\left[\begin{array}{cc}
                I_n & 0 \\
                0   & K
            \end{array}\right].
    $$
    Gillis, Mehrmann and Sharma\cite{gms} report finding a real Hurwitz stable pencil at a squared distance of $2.46$. Starting from a randomly generated point, our algorithm finds a real Hurwitz stable pencil at the squared distance of $1.01$. This pencil is saved in the variables $S$ and $T$ of the \texttt{example6dot2.mat} file that can be found in the GitHub repository. \footnote{See the previous footnote.} The computed pencil has no infinite eigenvalues, and is therefore in the feasible region for the DH method. It has a Jordan block of size four corresponding to a zero eigenvalue.
\end{example}

% \subsection{Statistical experiments} \label{sec:statistical}

Next, we perform more rigorous statistical comparisons with the DH method presented by Gillis, Mehrmann and Sharma\cite{gms}. In their work, the authors propose two variants of the method, which they call the GM method and the FGM method. Of these, the FGM method performs better based on heuristic observations. As such, we compare our method against the FGM method. 

We use the identity as a starting point to our method and optimize over $O(n) \times O(n)$, which guarantees a real output. For sizes $4 \leq n \leq 100$, we generate $10^2$ real Gaussian square pencils $A + \lambda B$ such that each element of $A$ and $B$ is drawn from the normal distribution $\mathcal{N}(0,1)$. We scale the pencil by a factor $(\sqrt{2} n)^{-1}$ in order to normalize the expectation of the Frobenius norm of the pencil. Figure \ref{fig:DH_comparison} shows the average distance computed by both algorithms. %as well as the median running time.

% For \cite{gms}, we used the same options that were given in the example scripts, i.e.

% \texttt{options.maxiter = 1e6; 
% options.timemax = 2; 
% options.posdef = 1e-6;}.

For the FGM method, we used the same options that were given in the example script \texttt{test\_case2\_randn.m} that was used to generate the randomized experiment of Gillis, Mehrmann and Sharma's work\cite{gms}, i.e. \texttt{options.maxiter = 1e6;
    options.timemax = 10;
    options.posdef = 1e-6}. In the script \texttt{test\_case2\_randn.m}, the parameter \texttt{options.timemax} was set to 2, but a comment noted that the value of 10 was used in the paper. We opted for the same choice of the parameter as in the paper.

% For our method, the maximum amount of time was set to 5 seconds. The rest of the parameters were the default parameters of the algorithm.

To facilitate comparison, the maximum amount of time for our method was set to be the running time of the FGM method for each pencil. This guarantees that our algorithm will never use more time than the FGM method. The rest of the parameters were the default parameters of the algorithm. The results are shown in Figure \ref{fig:DH_comparison}. %The running time of the DH method shows irregular behaviour. Nonetheless, 
The method of this paper significantly outperforms the FGM method for all sizes of the input pencil. While the FGM method by design does not find pencils that have index higher than one and our method sometimes finds pencils with higher indices, this cannot alone explain the large discrepancy in the computed distance because (by Lemma \ref{lem:samedistance}) given a stable pencil with index higher than one it is always possible to find another stable pencil, with index at most one, and at an arbitrarily small distance.

\begin{figure}[h!]
    \centering
    \begin{tikzpicture}
        \begin{axis}[
                legend pos = north west,
                width=\linewidth,
                yticklabel style={
                    /pgf/number format/fixed
                    % /pgf/number format/precision=5
                },
                scaled y ticks=false,
               height= 6cm,
                xlabel={rank of $B$},
                ylabel={distance}
            ]
             \addplot[red, dashed] table[y index = 1] 
            {statistical_comp_r.dat};
             \addplot[blue, dashed] table[y index = 2] 
            {statistical_comp_r.dat};
            \addplot[only marks, red, mark size =1pt] table[y index = 1] 
            {statistical_comp_r.dat};
            \addplot[only marks, blue, dashed, mark size =1pt] table[y index = 2] 
              {statistical_comp_r.dat};
            \legend{{\small{Riemann}}, {\small{FGM}}}
          \end{axis}
        \end{tikzpicture}

    % \begin{tikzpicture}
    % \begin{axis}[
    %         legend pos = south east,
    %         width=\linewidth,
    %         % xtick={5,10, 15, 20, 25, 30},
    %         %xticklabels={},
    %       % ymax=100,
    %        height= 6cm,
    %        title={{\small{Median running time}}},
    %         xlabel={rank of $B$},
    %         ylabel={time (s)}
    %     ]
    %      \addplot[red, dashed] table[y index = 3] 
    %     {statistical_comp_r.dat};
    %      \addplot[blue, dashed] table[y index = 4] 
    %     {statistical_comp_r.dat};
    %     \addplot[only marks, red, mark size =1pt] table[y index = 3] 
    %     {statistical_comp_r.dat};
    %     \addplot[only marks, blue, dashed, mark size =1pt] table[y index = 4] 
    %       {statistical_comp_r.dat};
    %     \legend{{\small{Riemann}}, {\small{DH}}}
    %   \end{axis}
    % \end{tikzpicture}

    \caption{Statistical comparison between the approach of this paper (labelled as Riemann in the picture) and the FGM algorithm \cite{gms} for $1 \leq r \leq 19$.} %\textcolor{red}{We can try to increase the maxtime (since it always gets triggered), but 10 seconds is the default for maxtime for their algorithm.}}
    \label{fig:DH_comparison_r}
\end{figure}
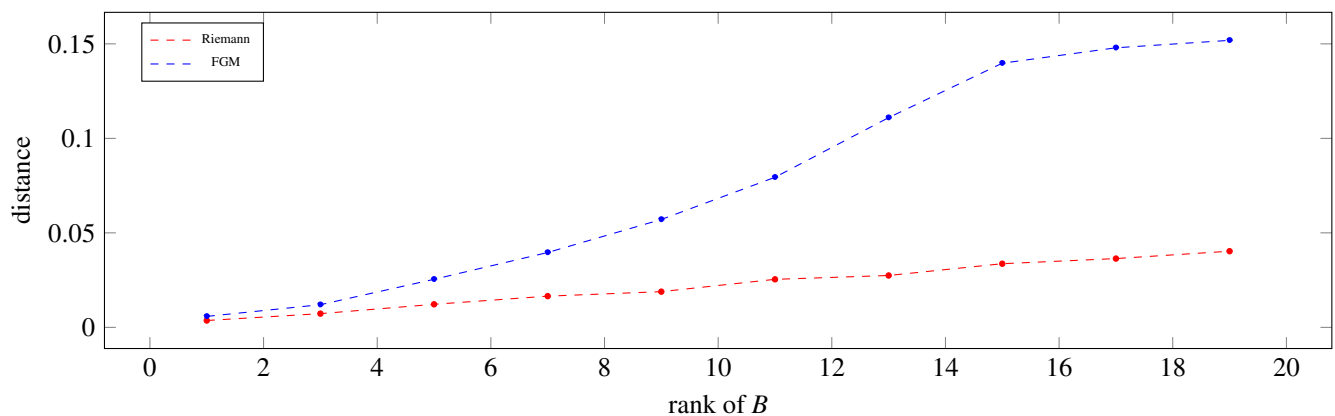

Finally, we compare the method of this paper with the FGM algorithm for pencils $A+xB$ where $B$ is of prescribed rank $r$. To do this, we randomly generate $10^2$ pencils $A+xB$ of size $n=20$ the same way as before and for each generated pencil we perform a rank-$r$ approximation for $B$ via SVD. We compare the methods for these input, for ranks $1 \leq r \leq 19$. The results are shown in Figure \ref{fig:DH_comparison_r}. Again, the method of this paper outperforms the FGM method significantly for all ranks $r$ of the matrix $B$. %This time, the median running time of the DH method is the maximum 10 seconds for all ranks $r$.

% \subsection{Gradient}
% To compute the gradient, we only need to modify the contribution coming from $g$:
% \begin{align*}
% \tilde g((a,b)) = 
% \begin{cases}
% \frac{(\|a\|-\|b\|)^2}{2} & \|a\| > \|b\|\\
% 0 & \|a\| \leq \|b\|
% \end{cases}.
% \end{align*}

% When $\|a\| > \|b\|$,
% \begin{align*}
%     D \tilde g((a,b))[(c,d)] &= (\|a\|-\|b\|) \text{Re} ([\frac{a^*}{\|a\|},-\frac{b^*}{\|b\|}] [c \ d]^T) 
% \end{align*}

% The directional derivative of $d_i(Q) := Q_{i,\cdot} (A+ \lambda B) Z_{\cdot,i}$ is 
% \begin{align*}
%     D d_i(Q)[X] = X_{i,\cdot} (A+ x B) Z_{\cdot,i}.
% \end{align*}
% The directional derivative of $Q \mapsto \tilde g(d_i(Q))$ is 
% \begin{small}
% \begin{align*}
%     &D(Q \mapsto \tilde g(d_i(Q)))(Q)[X] = (D((a,b) \mapsto \tilde g((a,b)))(d_i(Q)) \circ D d_i(Q))[X] \\
%     &= (\|(d_i)_a\|-\|(d_i)_b\|) \text{Re} ([\frac{(d_i)_a^*}{\|(d_i)_a\|},-\frac{(d_i)_b^*}{\|(d_i)_b\|}] [X_{i,\cdot} A Z_{\cdot,i}, X_{i,\cdot} B Z_{\cdot,i}]^T) \\
%     &=  (\|(d_i)_a\|-\|(d_i)_b\|) \text{Re} (X_{i,\cdot} \frac{(d_i)_a^*}{\|(d_i)_a\|}A Z_{\cdot,i} - X_{i,\cdot} \frac{(d_i)_b^*}{\|(d_i)_b\|}B Z_{\cdot,i}).
% \end{align*}
% \end{small}

% The gradient with respect to the $i$th row of $Q$ is then 
% \begin{small}
% \begin{align*}
% \nabla \bar g_{Q_{i,\cdot}} ( Q_{i,\cdot} & (A+ x B) Z_{\cdot,i}) = \\
%     & (\|(d_i)_a\|-\|(d_i)_b\|) \left( \frac{(d_i)_a}{\|(d_i)_a\|} (A Z_{\cdot,i})^* - \frac{(d_i)_b}{\|(d_i)_b\|} (B Z_{\cdot,i})^* \right). 
% \end{align*}    
% \end{small}

% The gradient with respect to $Z$ can be computed similarly. 

\subsection{Schur stability}\label{sec:numschur}
In this case, to the best of our knowledge, there are no available algorithms to compare with. On the other hand, the general remarks for the Hurwitz case hold also for the Schur case. Namely:
\begin{enumerate}
    \item The output usually has a long Jordan chain. Again, these results are analogous to those of the matrix-variant of the problem presented in \cite{np}.
    \item For randomly generated (real or complex) 10 by 10 input, a Schur stable pencil is often found at around the same distance as a Hurwitz stable pencil, namely, at a distance of approximately 5\% of the norm of the input.
    \item The running time is similar between the Schur and Hurwitz cases.
\end{enumerate}

For reference, we next run the algorithm for finding the nearest Schur stable pencil for the pencils in Examples \ref{ex1} and \ref{ex2}.   

\begin{example} \label{ex3}
    We consider the same pencil as in Example \ref{ex1}. Starting from a randomly generated point, our algorithm finds a Schur stable pencil at the squared distance of $1.85$. This pencil is saved in the variables $S$ and $T$ of the \texttt{example6dot3.mat} file that can be found in the GitHub repository. The starting point is saved in the $x0$ variable.
\end{example}

\begin{example} \label{ex4}
We consider the same pencil as in Example \ref{ex2}. Starting from a randomly generated point, our algorithm finds a Schur stable pencil at the squared distance of $1.02$. This pencil is saved in the variables $S$ and $T$ of the \texttt{example6dot4.mat} file that can be found in the GitHub repository. The starting point is saved in the $x0$ variable.
\end{example}

\section{Conclusion}\label{sec:conclusions}

In this paper, we have introduced a Riemannian optimization method for finding the nearest stable pencil. This method is a generalization of the work done by Noferini and Poloni\cite{np} who addressed the matrix-variant of the problem. Our reformulation of the problem reduces it to an unconstrained minimization problem of an almost-everywhere smooth cost function over the manifold $U(n) \times U(n)$. Numerical experiments demonstrate the effectiveness of this approach, even for finding the nearest real stable pencil, a variant of the problem for which our method was not specifically designed for. Future directions of research include developing a method specifically for the real-variant of the problem. 

\section*{Acknowledgements}

 We thank Federico Poloni for many useful discussions.

\end{document}